\documentclass {article}
\usepackage{authblk}
\usepackage{etex}
\usepackage[utf8]{inputenc}
\usepackage[english]{babel}
\usepackage{amsmath}
\usepackage{amsthm}
\usepackage{amssymb}
\usepackage{sectsty}
\usepackage{titlesec}
\usepackage{color}
\usepackage[color,matrix,arrow]{xy}
\usepackage{amsgen}
\usepackage{amstext}
\usepackage{amsbsy}
\usepackage{amsopn}
\usepackage{amsfonts}
\usepackage{eepic}
\usepackage{graphicx}
\usepackage{epsf}
\usepackage{pstricks}
\usepackage{enumerate}
\usepackage{eqnarray}
\usepackage{faktor}
\xyoption{all}

\usepackage{pgf}
\usepackage{tikz-cd}
\usetikzlibrary{automata, arrows.meta, positioning,decorations.pathmorphing}

\newcommand{\footrecall}[1]{%
} 
\usetikzlibrary{snakes,shapes,arrows,automata}

\newcommand\numberthis{\addtocounter{equation}{1}\tag{\theequation}}

\titleformat*{\section}{\large\bfseries}
\titleformat*{\subsection}{\normalsize \bfseries}

\newcommand{\N}{\mathbb{N}}

\newcommand{\Cyc}{\text{Cyc}}
\newcommand{\Rat}{\text{Rat}}

\newcommand{\Rec}{\text{Rec}}
\newcommand{\Alg}{\text{Alg}}
\newcommand{\CF}{\text{CF}}

\newcommand{\mc}{\mathcal}
\def\oo{\overline}

\theoremstyle{definition}
\newtheorem{theorem}{Theorem}[section]
\newtheorem{corollary}[theorem]{Corollary}

\newtheorem{proposition}[theorem]{Proposition}

\newtheorem{lemma}[theorem]{Lemma}

\newtheorem{remark}[theorem]{Remark}

\begin{document}

\title{Geodesic languages for rational subsets and conjugates in virtually free groups}
\author{Andr\'e Carvalho\thanks{andrecruzcarvalho@gmail.com}\, and Pedro V. Silva\thanks{pvsilva@fc.up.pt}}
\affil{Centre of Mathematics, Department of Mathematics

Faculty of Sciences of the University of Porto

 R. Campo Alegre, 4169-007 Porto,
 
Portugal}
\maketitle

\begin{abstract}
We prove that a subset of a virtually free group is rational if and only if the language of geodesic words representing its elements (in any generating set) is rational and that the language of geodesics representing conjugates of elements in a rational subset of a virtually free group is context-free. As a corollary, the doubly generalized conjugacy problem is decidable for rational subsets of finitely generated virtually free groups: there is an algorithm taking as input two rational  subsets $K_1$ and $K_2$ of a virtually free group that decides whether there is one element of $K_1$ conjugate to an element of $K_2$. For free groups, we prove that the same problem is decidable with rational constraints on the set of conjugators.
\end{abstract}

\section{Introduction}
Given a group $G$, two elements $x,y\in G$ are said to be \emph{conjugate} if there is some $z\in G$ such that $x=z^{-1}yz$, in which case we write $x\sim y$. The \emph{conjugacy problem} CP$(G)$ consists on, given $x,y\in G$, deciding whether $x\sim y$ or not. This was one of the three algorithmic problems introduced by Dehn \cite{[Deh11]}, together with the \emph{word problem} and the \emph{isomorphism problem}.

The word problem, WP$(G)$, is possibly the most well-studied algorithmic problem in group theory and consists on, given a word on the generators of a group, deciding whether the element represented by that word is the identity or not, or, equivalently, given two words on the generators, deciding whether they represent the same group element. The \emph{membership problem}, MP$(G)$, also known as the \emph{generalized word problem} consists on, given a finitely generated subgroup $H\leq G$ and an element $x\in G$, deciding whether $x\in H$ or not. This can be considered more generally for subsets belonging to a reasonably well-behaved class instead of subgroups (e.g. rational or context-free subsets). This can also be rewritten as the question of deciding whether there is some $y\in H$ such that $x=y$ (see \cite{[Loh15]} for a survey on this problem). In the same spirit, a generalization of the conjugacy problem was considered in \cite{[LS11]} and proven to be decidable with respect to rational subsets of finitely generated virtually free groups. The \emph{generalized conjugacy problem with respect to $\mc C$}, GCP$_\mc C(G)$, where $\mc C$ is a class of subsets of $G$ consists then on, given $x\in G$ and $K\in \mc C$, deciding whether there is some $y\in K$ such that $x\sim y$. Clearly, if $\mc C$ contains all singletons (which occurs if $\mc C$ is the class of rational subsets or the class of cosets of finitely generated subgroups), this is indeed a generalization of the conjugacy problem.

The \emph{intersection problem} IP$_\mc C(G)$ consists on, given two subsets $K_1,K_2\in \mc C$, deciding whether $K_1\cap K_2=\emptyset$. Naturally, if $\mc C$ is a class of subsets containing all singletons, if we can decide the intersection problem with respect to $\mc C$, we can decide the membership problem with respect to $\mc C$. Thus, in some sense, the intersection problem can be seen as the doubly generalized word problem and, as done above, it can also be rewritten as the question of deciding  whether there are some $x\in K_2$ and $y\in K_2$ such that $x=y$.  However, if the class of subsets is closed under product of subsets and inversion, this is equivalent to the membership problem, as it consists on deciding whether $1\in K_1K_2^{-1}$. This is the case when considering rational or algebraic subsets, but does not hold in general. In this paper, we consider the \emph{doubly generalized conjugacy problem with respect to $\mc C$}, DGCP$_{\mc C}(G)$, which is the natural generalization of the conjugacy problem corresponding to the intersection problem, that is, the problem of, given $K_1,K_2\in \mc C$,  deciding whether there are some $x\in K_2$ and $y\in K_2$ such that $x\sim y$.

In case $\mc C$ is the class of the rational subsets of $G$, the following is easy to see (where $\leq$ means that the problem on the left-hand side is reducible to the one on the right-hand side and $\equiv$ means that the problems are equivalent): 
\begin{align*}
&\text{WP}(G)&&\leq &&&\text{MP}_{\Rat}(G)&&&&\equiv &&&&&\text{IP}_{\Rat}(G)\\
&\quad\rotatebox{270}{$\leq\;\;$}&&  &&&\rotatebox{270}{$\leq\;\;$}\quad\;\;&&&& &&&&&\quad\rotatebox{270}{$\leq\;\;$}\\
&\text{CP}(G)&&\leq &&&\text{GCP}_{\Rat}(G)&&&&\leq &&&&&\text{DGCP}_{\Rat}(G)
\end{align*}

We will additionally consider versions of the conjugacy problems with certain \emph{constraints} on the conjugators. In \cite{[LS11]}, it is proved that the \emph{generalized conjugacy problem with rational constraints} with respect to rational subsets of finitely generated virtually free groups is decidable, meaning that, given a virtually free group $G$, there is an algorithm taking as input two rational subsets $L,K\in \Rat(G)$ and an element $x\in G$ and decides if there is some $z\in L$ such that $z^{-1}xz\in K$.

Given $K,L \subseteq G$, let 
$$\alpha(K,L) = \bigcup_{u \in L} u^{-1}Ku.$$
When $L=G$, we simply write $\alpha(K)$ to denote  $\alpha(K,G)$ 

In this paper, we will present a language-theoretical proof of the decidability of the doubly generalized conjugacy problem with rational constraints with respect to rational subsets of finitely generated free groups. To do so, we prove that, in a finitely generated free group,  the set $\alpha(K,L)$ of all elements conjugate to an element of $K$ by an element of $L$ is  a context-free subset of the ambient free group. 
It is proved in \cite{[Lev23]} that a group is virtually free if and only if its conjugacy classes are context-free subsets. Equivalently, a group is virtually free if and only if $\alpha(S)$ is context-free for all singletons $S$. We prove something much stronger in the case of free groups, namely that $\alpha(K,L)$ is context-free if both $K$ and $L$ are rational.

\newtheorem*{main}{Theorem \ref{main}}
\begin{main}
Let $K,L \in {\rm Rat}\,F_A$. Then $\alpha(K,L)$ is a context-free subset of $F_A$.
\end{main}

Since context-free languages are closed under intersection with regular languages and emptiness of context-free languages is decidable, we have the following corollary:
\newtheorem*{dgcprat free}{Corollary \ref{dgcprat free}}
\begin{dgcprat free}
The doubly generalized conjugacy problem with rational constraints  is decidable with respect to rational subsets of a finitely generated free group.
\end{dgcprat free}

Regarding virtually free groups, we prove a generalization of the well-known Benois's theorem, showing that a subset is rational if and only if the language of geodesics representing its elements is rational.

\newtheorem*{benois vfree}{Corollary \ref{benois vfree}}
\begin{benois vfree}
Let $G$ be a f.g. virtually free group and $K\subseteq G$. The following are equivalent:
\begin{enumerate}
\item $K\in \Rat(G)$.
\item $Geo_X(K)$ is a rational language for some finite generating set $X$.
\item $Geo_X(K)$ is a rational language for every finite generating set $X$.
\end{enumerate} 
Moreover, the constructions are effective.
\end{benois vfree}

We then prove that the language of geodesics representing conjugates of a given rational subset is context-free, which yields  a language-theoretic proof of the decidability of the  doubly generalized conjugacy problem for rational subsets of finitely generated virtually free groups.
Again, it follows from \cite{[Lev23]} that, for a virtually free group and a singleton $S\subseteq G$, $Geo(\alpha(S))=S\pi^{-1}\cap Geo(G)$ is a context-free language. We prove that this holds for every rational subset of $G$. 
\newtheorem*{alpha cf}{Theorem \ref{alpha cf}}
\begin{alpha cf}
Let $G$ be a virtually free group and $K\in \Rat(G)$. Then $Geo(\alpha(K))$ is context-free.
\end{alpha cf}

\newtheorem*{dgcp vfree}{Corollary \ref{dgcp vfree}}
\begin{dgcp vfree}
Let $G$ be a virtually free group. Then the doubly generalized conjugacy problem is decidable.
\end{dgcp vfree}

We remark that Corollaries \ref{dgcprat free} and \ref{dgcp vfree} were already known, as it follows directly from the fact that the existential theory of equations with rational constraints in free groups is PSPACE-complete, which was proved by Diekert, Guti\'errez and Hagenah in \cite{[DGH05]}, since, for $L,K_1,K_2\in \Rat(G)$, the statement that there is an element of $K_1$ conjugate to an element of $K_2$ by an element of $L$ can be expressed as: 
$$\exists z\in L\, \exists x\in K_1\, \exists y\in K_2 \,:\, z^{-1}xz=y.$$
Similarly, for virtually free groups, we can use the analogous result for virtually free groups proved in \cite{[DG10]}.

\section{Preliminaries}\label{sec:prelim}
In this section, we will present basic definitions and results on rational, algebraic and context-free subsets of groups (for more details, the reader is referred to \cite{[Ber79]} and \cite{[BS21]}) and on virtually free groups.

\subsection{Subsets of groups}

The set $\{1,\ldots, n\}$ will be denoted by $[n]$.
Let $G=\langle A\rangle$ be a finitely generated group, $A$ be a finite generating set, $\widetilde A=A\cup A^{-1}$ and $\pi:\widetilde A^*\to G$ be the canonical (surjective) homomorphism. This notation will be kept throughout the paper.

A subset $K\subseteq G$ is said to be \emph{rational} if there is some rational language $L\subseteq \widetilde A^*$ such that $L\pi=K$ and \emph{recognizable} if $K\pi^{-1}$ is rational. 

We will denote by $\Rat(G)$ and $\Rec(G)$ the class of rational and recognizable subsets of $G$, respectively. Rational subsets generalize the notion of finitely generated subgroups.

\begin{theorem}[\cite{[Ber79]}, Theorem III.2.7]
\label{AnisimovSeifert}
Let $H$ be a subgroup of a group $G$. Then $H\in \Rat(G)$ if and only if $H$ is finitely generated.
\end{theorem}

Similarly, recognizable subsets generalize the notion of finite index subgroups.

\begin{proposition}
\label{rec fi}
Let $H$ be a subgroup of a group $G$. Then $H\in \Rec(G)$ if and only if $H$ has finite index in $G$.
\end{proposition}

In fact, if $G$ is a group and $K$ is a subset of $G$ then $K$ is recognizable if and only if $K$ is a (finite) union of cosets of a subgroup of finite index.

In case the group $G$ is a free group with basis $A$ with surjective homomorphism $\pi:\widetilde A^*\to G$, given  $L\subseteq \widetilde A^*$, we define the set of reduced words representing elements in $L\pi$ by $$\oo L=\{w\in \widetilde A^*\mid w \text{ is reduced and there exists $u\in L$ such that $u\pi=w\pi$} \}.$$ 
Benois' Theorem provides us with a useful characterization of rational subsets in terms of reduced words representing the elements in the subset.
\begin{theorem}[Benois]
Let $F$ be a finitely generated free group with basis $A$ and let $L \subseteq \widetilde{A}^*$. Then $\oo{L}$ is a rational language of $\widetilde{A}^*$ if and only if $L\pi$ is a rational subset of $F$. 
\end{theorem}

A natural generalization of these concepts concerns the class of context-free languages. 
A subset $K\subseteq G$ is said to be \emph{algebraic} if there is some context-free language $L\subseteq \widetilde A^*$ such that $L\pi=K$ and \emph{context-free} if $K\pi^{-1}$ is context-free. 
We will denote by $\Alg(G)$ and $\CF(G)$ the class of algebraic and context-free subsets of $G$, respectively.  It follows from \cite[Lemma 2.1]{[Her91]} that these definitions, as well as the definitions of rational and recognizable subsets, do not depend on the finite alphabet $A$ or the surjective homomorphism $\pi$.

It is obvious from the definitions that $\Rec(G)$, $\Rat(G)$, $\CF(G)$ and $\Alg(G)$ are closed under union, since both rational and context-free languages are closed under union. The intersection case is distinct: from the fact that rational languages are closed under intersection, it follows that $\Rec(G)$ must be closed under intersection too. However $\Rat(G)$, $\Alg(G)$ and $\CF(G)$ might not be. Another important closure property is given by the following lemma from \cite{[Her91]}.

\begin{lemma}\cite[Lemma 4.1]{[Her91]}
\label{lemma herbst}
 Let $G$ be a finitely generated group, $R\in \Rat(G)$ and $C \in \{\Rec, \CF\}$.
If $K \in C(G)$, then $KR, RK \in C(G).$
\end{lemma}
The following is an immediate consequence of the previous lemma.
\begin{corollary}
\label{conj}
Let $K \subseteq F_A$ and $u \in F_A$. Then $u^{-1}Ku$ is context-free if and only if $K$ is context-free.
\end{corollary}

For a finitely generated group $G$, it is immediate from the definitions that $$\Rec(G)\subseteq \CF(G) \subseteq \Alg(G)$$ and that $$\Rec(G)\subseteq \Rat(G) \subseteq \Alg(G).$$

\begin{lemma}\cite[Lemma 4.3]{[Her91]}\label{square}
Let $X,Y$ be finite alphabets and let $\psi:Y^*\to M$, $\varphi: X^*\to M'$ be homomorphisms onto monoids $M,M'$. Then every homomorphism $\tau:M'\to M$ can be lifted to a homomorphism $h:X^*\to Y^*$ such that the diagram
$$
\begin{tikzcd}[sep=large]
X^*     \ar[r,dashed,"h"] \ar[d,swap,"\varphi"] & Y^*            \ar[d,"\psi"]                \\
 M'\ar[r,swap,"\tau"] & M
\end{tikzcd}
$$
commutes. As a consequence, $T\tau^{-1}\varphi^{-1} = T\psi^{-1}h^{-1}$ for every $T \subseteq M$.
\end{lemma}

However, there is no general inclusion between $\Rat(G)$ and $\CF(G)$. For example, if $G$ is virtually abelian, then $\CF(G)\subseteq \Alg(G)= \Rat(G)$ (and the inclusion is strict if the group is not virtually cyclic) and if the group is virtually free, then $\Rat(G)\subseteq \CF(G)$ (see \cite[Lemma 4.2]{[Her91]}).

In the case of the free group $F_n$ of rank $n\geq 1$, Herbst proves in  \cite{[Her91]} an analogue of Benois's Theorem for context-free subsets:

\begin{lemma} \cite[Lemma 4.6]{[Her91]}  \label{herbst benois cf}
Let $F$ be a finitely generated free group with basis $A$ and let $L \subseteq \widetilde{A}^*$. Then $\oo{L}$ is a context-free language of $\widetilde{A}^*$ if and only if $L\pi$ is a context-free subset of $F$. 
\end{lemma}

A slight improvement of the previous lemma can be easily obtained:

\begin{lemma} 
\label{generalized benois cf}
Let $F$ be a finitely generated free group and $K\subseteq F$. Then $K\in \CF(F)$ if and only if there is a context-free language $\oo K\subseteq L\subseteq K\pi^{-1}$.
\end{lemma}

We will also make use of the following lemma, which is a simple exercise:

\begin{lemma}
\label{cfl}
Let $L \in {\rm Rat} \, A^*$ and let $u,v \in A^*$. Then the languages
$$\bigcup_{n \geq 0} u^nLv^n, \quad \bigcup_{0 \leq m \leq n} u^mLv^n \quad \mbox{and} \quad \bigcup_{0 \leq m \leq n} u^nLv^m$$
are all context-free.
\end{lemma}

\subsection{Virtually free groups}
A group $G$ is said to be \emph{virtually free} if it has a free subgroup $F$ of finite index. Since subgroups of free groups are free and every finite index subgroup contains a finite index normal subgroup, we can assume that $F\trianglelefteq_{f.i.} G$. We will usually write $$G=Fb_1\cup\cdots \cup Fb_n,$$
where all cosets $Fb_i$ are disjoint.

Algebraic and context-free subsets of virtually free groups are studied in \cite{[Car23c]}. In particular, it is proved in \cite[Theorem 4.3]{[Car23c]} that, if $G$ is a finitely generated virtually free group and $H\leq_{f.g.} G$, then
\begin{align} \label{fatou vfree}
\CF(H)=\{K\subseteq H \mid K\in \CF(G)\}.
\end{align}

Also, combining \cite[Proposition 4.1]{[Sil02b]} and Propositions 3.6 and 3.7 of \cite{[Car23c]} we have that $\Rec(G)$ (resp. $\Rat(G)$, $\Alg(G)$, $\CF(G)$) consists of sets of the form 
$L_ib_i$, where $L_i\in \Rec(F)$ (resp. $\Rat(F)$, $\Alg(F)$, $\CF(F)$).

A word $u=u_1\ldots u_k$ is said to be \emph{cyclically reduced} if $u_1\neq u_k^{-1}$. Every word $u$ can be written as $u=w^{-1}\widetilde uw$, where $\widetilde u$ is cyclically reduced. We refer to $\widetilde u$ as the \emph{cyclically reduced core} of $u$.

 \section{Conjugates of elements in a rational subset}\label{sec: DGCP}

 In this section, we will  prove that, in a free group, the set of conjugates of elements in a rational subset $K$ with a conjugator in a rational subset $L$, $\alpha(K,L)$, is context free and that a context-free grammar representing it can be effectively computed. As a corollary, we have 
   a language-theoretical proof that the doubly generalized conjugacy problem with respect to rational subsets with rational constraints is decidable on a free group $F$. This result also follows from the very strong theorem by Diekert, Gutiérrez and Hagenah \cite{[DGH05]} stating that the existential theory of equations with rational constraints in free groups is PSPACE-complete. 

 We will then consider the case of virtually free groups. We start by proving a generalization of Benois's theorem: a subset of a virtually free group is rational if and only if the language of geodesic words representing its elements is rational. Then, we show that the language of geodesic words representing a conjugate of an element in a rational subset $K$ is context-free (and computable), obtaining a language-theoretical proof for the doubly generalized conjugacy problem. This problem was already known to be decidable, as its decidability follows from the solution of equations with rational constraints for virtually free groups \cite{[DG10]}. 

\subsection{Free groups}

Given $K,L \subseteq F_A$, we say that the product $KL$ is reduced if $\oo{K}\,\oo{L} \subseteq \oo{\widetilde{A}^*}$. The purpose of this subsection is to prove that, 
given $K,L \in {\rm Rat}\,F_A$, then $\alpha(K,L)$ is a context-free subset of $F_A$. We start by solving the particular cases where $L$ and $K$ satisfy some reducibility conditions.

\begin{lemma}
\label{lang}
Let $K,L \in {\rm Rat}\, \widetilde{A}^*$. Then $\bigcup_{u \in L} u^{-1}Ku \subseteq \widetilde{A}^*$ is a context-free language.
\end{lemma}

\noindent\textit{Proof.}
Let ${\cal{G}} = (V,P,S)$ be the context-free grammar on the alphabet $\widetilde{A} \cup \{ \$ \}$ defined by $V = \widetilde{A} \cup \{ \$, S \}$ and $P = \{ (S,aSa^{-1}) \mid a \in \widetilde{A}\} \cup \{ (S,\$)\}$. It is immediate that $L({\cal{G}}) = \{ u^{-1}\$ u \mid u \in \widetilde{A}^* \}$. Since context-free languages are closed under intersection with regular languages, it follows that 
$$\{ u^{-1}\$ u \mid u \in L \} = \{ u^{-1}\$ u \mid u \in \widetilde{A}^* \} \cap L^{-1}\$ L$$
is context-free. Since context-free languages are closed under substitution, we can replace the letter $\$ $ by the rational (hence context-free) language $K$ and remain context-free. Therefore $\bigcup_{u \in L} u^{-1}Ku$ is a context-free language.
\qed\\

\begin{lemma}
\label{red}
Let $K,L \in {\rm Rat}\,F_A$ with both $L^{-1}\, K$ and $KL$ reduced. Then $\alpha(K,L)$ is a context-free subset of $F_A$.
\end{lemma}

\noindent\textit{Proof.}
By Lemma \ref{herbst benois cf}, it suffices to show that 
$\oo{\alpha(K,L)}$ is a context-free language. This same argument will be used in the next proofs without further reference.

Now $\oo{\alpha(K,L)} = \{ u^{-1}\oo{K}u \mid u \in \oo{L} \}$ and it follows from Benois' Theorem 
that $\oo{K}$ and $\oo{L}$ are both rational languages. 
 By Lemma \ref{lang}, $\oo{\alpha(K,L)}$ is a context-free language.
\qed\\

\begin{lemma}
\label{rred}
Let $K,L \in {\rm Rat}\,F_A$ with $KL$ reduced. Then $\alpha(K,L)$ is a context-free subset of $F_A$.
\end{lemma}
\noindent\textit{Proof.}
We may assume that $K$ and $L$ are both nonempty.
Let $C$ denote the set of all cyclically reduced elements of $F_A$, which is clearly a rational subset. Then $K \cap C$ and $K \setminus C$ are both rational subsets of $F_A$. Since $L^{-1}(K\setminus C)$ and $(K\setminus C)L$ are both reduced, it follows from Lemma \ref{red} that $\alpha(K\setminus C,L)$ is a context-free subset of $F_A$. Thus it suffices to show that $\alpha(K\cap C,L)$ is a context-free subset of $F_A$. Therefore we may assume that $K \subseteq C$, and we may also assume that $1 \notin K$.

Let ${\cal{A}} = (Q,q_0,T,E)$ and ${\cal{A}}' = (Q',q'_0,T',E')$ denote respectively the minimal automata of $\oo{L}$ and $\oo{K}$. For all $I,J \subseteq Q$, let 
$L_{IJ} = L(Q,I,J,E)$. For all $I',J' \subseteq Q'$, let 
$L'_{I'J'} = L(Q',I',J',E')$. Let
$$X = \bigcup_{m = 0}^{|Q|-1} Q^{2m+1} \times Q'.$$
We show that
\begin{equation}
\label{rred1}
\begin{array}{lll}
\oo{\alpha(K,L)}&=&\{ w^{-1}v_2v_1w \mid \exists\,(p_1,q_1,\ldots,p_m,q_m,p_{m+1},q') \in X,\, v_1 \in L'_{q'_0q'} \cap (\bigcap_{i=0}^m L_{q_ip_{i+1}}),\\ &&\\
&&v_2 \in L'_{q'T'} \cap (\bigcap_{i=1}^m L_{p_iq_{i}})
,\, w \in L_{p_{m+1}T} \} \cap \oo{\widetilde{A}^*}.
\end{array}  
\end{equation}

Indeed, let $w^{-1}v_2v_1w$ belong to the right hand side of (\ref{rred1}) for some $(p_1,q_1,\ldots,p_m,q_m,p_{m+1},q') \in X$. Then we have a path $q'_0 \xrightarrow{v_1} q' \xrightarrow{v_2} t' \in T'$ in $\cal{A}'$ and a path
$$q_0 \xrightarrow{v_1} p_1 \xrightarrow{v_2} q_1 \xrightarrow{v_1} \ldots \xrightarrow{v_2} q_{m-1} \xrightarrow{v_1}  p_m \xrightarrow{v_2} q_m \xrightarrow{v_1} p_{m+1} \xrightarrow{w} t \in T$$
in $\cal{A}$.  Hence $v_1v_2 \in L({\cal{A}'}) = \oo{K}$ and $(v_1v_2)^mv_1w \in L({\cal{A}}) = \oo{L}$. It follows that
$$w^{-1}v_2v_1w = w^{-1}v_1^{-1}(v_2^{-1}v_1^{-1})^mv_1v_2(v_1v_2)^mv_1w \in \alpha(K,L).$$
Since $w^{-1}v_2v_1w$ is a reduced word by hypothesis, we get $w^{-1}v_2v_1w \in \oo{\alpha(K,L)}$.

Conversely, assume that $u \in \oo{L}$ and $v \in \oo{K}$. We consider the longest prefix of $u$ which is a prefix of some power of $v$. More precisely, write $u = v^mv_1w$ such that $m \geq 0$ and $v = v_1v_2$ with $v_2 \neq 1$. Then
$$\oo{u^{-1}vu} = \oo{w^{-1}v_1^{-1}v^{-m}vv^mv_1w} = \oo{w^{-1}v_2v_1w}.$$
Note that every path of the form $p \xrightarrow{v^m} q$ in $\cal{A}$ with $m \geq |Q|$ must contain some loop labelled by $v^s$ with $1 \leq s \leq |Q|$, hence we may replace $u = v^mv_1w$ by $u' = v^{m-s}v_1w$ without changing the final outcome $\oo{w^{-1}v_2v_1w}$. Thus we assume that $m < |Q|$. We must have a path
$$q_0 \xrightarrow{v_1} p_1 \xrightarrow{v_2} q_1 \xrightarrow{v_1} \ldots \xrightarrow{v_2} q_{m-1} \xrightarrow{v_1}  p_m \xrightarrow{v_2} q_m \xrightarrow{v_1} p_{m+1} \xrightarrow{w} t \in T$$
in $\cal{A}$ and a path $q'_0 \xrightarrow{v_1} q' \xrightarrow{v_2} t' \in T'$ in $\cal{A}'$. It follows that
$(p_1,q_1,\ldots,p_m,q_m,p_{m+1},q') \in X$, $v_1 \in L'_{q'_0q'} \cap (\bigcap_{i=0}^m L_{q_ip_{i+1}})$, $v_2 \in L'_{q'T'} \cap (\bigcap_{i=1}^m L_{p_iq_{i}})$ and $w \in L_{p_{m+1}T}$. It remains to show that $w^{-1}v_2v_1w$ is reduced.

Indeed, $v_2v_1w$  
labels a path in a trim automaton recognizing a reduced language, hence must be a reduced word itself. Suppose that 
$w^{-1}v_2$ is not reduced. Then $v_2$ and $w_2$ share the same first letter, say $a$. Then $v^mv_1a$ is a prefix of $u$ which is a prefix of $v^{m+1}$, contradicting the maximality of $v^mv_1$.
Hence $w^{-1} v_2$ is reduced. Since $v_2v_1w$ is reduced and $v_2 \neq 1$, then $w^{-1} v_2v_1w$ is itself reduced and so  (\ref{rred1}) holds. 
 Now, applying Lemma \ref{lang} $|X|$ times to the rational languages featuring the right hand side of (\ref{rred1}),
 and taking into account that context-free languages are closed under intersection with rational languages and union, we conclude that $\oo{\alpha(K,L)}$ is a context-free language as intended.
\qed\\

\begin{lemma}
\label{lred}
Let $K,L \in {\rm Rat}\,F_A$ with $L^{-1}K$ reduced. Then $\alpha(K,L)$ is a context-free subset of $F_A$.
\end{lemma}
\noindent\textit{Proof.}
 We may assume that $K$ and $L$ are both nonempty.
Since $K$ rational implies $K^{-1}$ rational and $L^{-1}K$ reduced implies $K^{-1}L$ reduced, it follows from the proof of Lemma \ref{rred} that $\oo{\alpha(K^{-1},L)}$ is a context-free language. Now
$$\oo{\alpha(K,L)} = \bigcup_{u \in L} \oo{u^{-1}Ku} = (\bigcup_{u \in L} \oo{u^{-1}K^{-1}u})^{-1} = (\oo{\alpha(K^{-1},L)})^{-1}.$$
Since context-free languages are closed under reversal and homomorphism, it follows easily that $(\oo{\alpha(K^{-1},L)})^{-1}$ is a context-free language. Thus $\oo{\alpha(K,L)}$ is a context-free language and we are done.
\qed\\

Now, we can prove the main result of this subsection.

\begin{theorem}
\label{main}
Let $K,L \in {\rm Rat}\,F_A$. Then $\alpha(K,L)$ is a context-free subset of $F_A$.
\end{theorem}
\noindent\textit{Proof.}
 We may assume that $K$ and $L$ are both nonempty.
Since $\alpha(1,L) = 1$, we may assume that $1 \notin K$.

Let ${\cal{A}} = (Q,q_0,T,E)$ and ${\cal{A}}' = (Q',q'_0,T',E')$ denote respectively the minimal automata of $\oo{L}$ and $\oo{K}$. We keep the notation introduced in the proof of Lemma \ref{rred}. 
We define
$$X = \{ (q,p',q') \in Q \times Q' \times Q' \mid L_{q_0q} \cap L'_{q'_0p'} \cap (L'_{q'T'})^{-1},\, L'_{p'q'} \setminus \{1\} \neq \emptyset\}.$$
For every $a \in \widetilde{A}$, we define  the possibly empty subsets of $F_A$
$$Y_a = \{ w^{-1}v_2w \mid \exists\,(q,p',q') \in X,\, v_2 \in L'_{p'q'} \cap \widetilde{A}^*a,\, w \in L_{qT} \setminus a^{-1}\widetilde{A}^*\},$$
$$Z_a = \{ w^{-1}v_2w \mid \exists\, (q,p',q') \in X,\, v_2 \in L'_{p'q'} \cap a\widetilde{A}^*,\, w \in L_{qT} \setminus a\widetilde{A}^*\}.$$
We show that
\begin{equation}
\label{main1}
\alpha(K,L) = \bigcup_{a \in \widetilde{A}} (Y_a \cup Z_a).
\end{equation}
Let $y \in Y_a$. Then there exist $(q,p',q') \in X$, $v_2 \in L'_{p'q'} \cap \widetilde{A}^*a$ and $w \in L_{qT} \setminus a^{-1}\widetilde{A}^*$ such that $y = w^{-1}v_2w$. Since $(q,p',q') \in X$, there exists some $v_1 \in L_{q_0q} \cap L'_{q'_0p'} \cap (L'_{q'T'})^{-1}$. Then we have a path $q_0 \xrightarrow{v_1} q \xrightarrow{w} t \in T$
in $\cal{A}$ and a path $$q'_0 \xrightarrow{v_1} p' \xrightarrow{v_2} q'  \xrightarrow{v_1^{-1}} t' \in T'$$ in $\cal{A}'$.
Hence $v_1v_2v_1^{-1} \in L({\cal{A}'}) = \oo{K}$ and $v_1w \in L({\cal{A}}) = \oo{L}$. It follows that
$$y = w^{-1}v_2w = w^{-1}v_1^{-1}(v_1v_2v_1^{-1})v_1w \in \alpha(K,L).$$
Thus $Y_a \subseteq \alpha(K,L)$. The inclusion $Z_a \subseteq \alpha(K,L)$ is proved similarly.

Conversely, assume that $u \in \oo{L}$ and $v \in \oo{K}$. We may write $v = bcb^{-1}$ with $c$ cyclically reduced. Since $1 \notin K$, we have $c \neq 1$. Let $v_1$ denote the longest common prefix of $u$ and $b$ 
 (which may be the empty word).
Write $u = v_1w$ and $v = v_1v_2v_1^{-1}$. We must have a path
$$q'_0 \xrightarrow{v_1} p' \xrightarrow{v_2} q' \xrightarrow{v_1^{-1}}  t' \in T'$$
in $\cal{A}'$ and a path $q_0 \xrightarrow{v_1} q \xrightarrow{w} t \in T$ in $\cal{A}$. It follows that
$v_1 \in L_{q_0q} \cap L'_{q'_0p'} \cap (L'_{q'T'})^{-1}$, $v_2 \in L'_{p'q'} \setminus \{ 1\}$ and $w \in L_{qT}$, hence $(q,p',q') \in X$.  Moreover, $u^{-1}vu = w^{-1}v_1^{-1}(v_1v_2v_1^{-1})v_1w = w^{-1}v_2w$.

Now it follows from the maximality of $v_1$ that at least one of the products 
$w^{-1}v_2, v_2w$ must be reduced
 (if $v_1 = b$, this follows from $v_2 = c$ being cyclically reduced).
If $w^{-1}v_2$ is reduced, then $u^{-1}vu \in Z_a$ when $a$ denotes the first letter of $v_2$. If $v_2w$ is reduced, then $u^{-1}vu \in Y_a$ when $a$ denotes the last letter of $v_2$. Therefore (\ref{main1}) holds.

Since
$$Y_a = \bigcup_{(q,p',q') \in X} \alpha(L'_{p'q'} \cap \widetilde{A}^*a, L_{qT} \setminus a^{-1}\widetilde{A}^*)$$
and $(L'_{p'q'} \cap \widetilde{A}^*a)(L_{qT} \setminus a^{-1}\widetilde{A}^*)$ is reduced, it follows from Lemma \ref{rred} that $Y_a$ is a context-free subset of $F_A$.

Since
$$Z_a = \bigcup_{(q,p',q') \in X} \alpha(L'_{p'q'} \cap a\widetilde{A}^*, L_{qT} \setminus a\widetilde{A}^*)$$
and $(L_{qT} \setminus a\widetilde{A}^*)^{-1}(L'_{p'q'} \cap a\widetilde{A}^*)$ is reduced, it follows from Lemma \ref{lred} that $Z_a$ is a context-free subset of $F_A$.

Now it follows from (\ref{main1}) that $\alpha(K,L)$ is itself a context-free subset of $F_A$.
\qed\\

The following corollary follows as an immediate application of the previous theorem and will be useful in the next subsection to deal with the case of virtually free groups.
\begin{corollary}
\label{powers}
Let $K \in {\rm Rat}\,F_A$ and $u \in F_A$. Then $\bigcup_{n\in \N} u^{-n}Ku^n$ is a context-free subset of $F_A$.
\end{corollary}
\noindent\textit{Proof.}
Since $u^* \in {\rm Rat}\,F_A$, the claim follows immediately from Theorem \ref{main}.
\qed\\

Since context-free languages are closed under intersection with regular languages and emptiness of a context-free language can be decided, we can decide the  doubly generalized conjugacy problem with rational constraints with respect to rational subsets of a finitely generated free groups.

\begin{corollary}\label{dgcprat free}
The doubly generalized conjugacy problem with rational constraints  is decidable with respect to rational subsets of a finitely generated free groups.
\end{corollary}
\noindent\textit{Proof.} Let $\pi:\widetilde A^*\to F_A$ be the canonical surjective homomorphism and $K_0, K_1,K_2\in \Rat(F_A)$ be our input (by this we mean that we get three finite state automata recognizing languages $L_0, L_1,L_2\subseteq \widetilde A^*$ such that $L_i\pi=K_i$, for $i=0,1,2.$ We want to decide if there are some $u\in K_0$, $x_1\in K_1$ $x_2\in K_2$ such that $x_1=u^{-1}x_2u$, i.e., if $K_1\cap \alpha(K_2,K_0)=\emptyset$. In view of Theorem \ref{main}, we can compute a context-free grammar $\mc G$ such that $L(\mc G)=(\alpha(K_2,K_0))\pi^{-1}$. Then $L_1\cap L(\mc G)$ is an effectively constructible context-free language and  $K_1\cap \alpha(K_2,K_0)=\emptyset$ if and only if  $ L_1\cap L(\mc G)=\emptyset$, which can be decided. 
\qed\\

\subsection{Virtually free groups}

Now we turn our attention to the case of virtually free groups. Our goal is to prove that the doubly generalized conjugacy problem is decidable with respect to rational subsets. We will write $G$ to denote a f.g. virtually free group and put $$G=Fb_1\cup\cdots\cup Fb_m,$$ where $F=F_A$ is a free normal  subgroup of $G$ of finite index $m$. We will also  put $B=A\cup \{b_1,\ldots, b_m\}$. Unless stated otherwise, $B$ will be our standard generating set for $G$.

For a subset $K\subseteq G$ and a generating set $X$ of $G$, let $Geo_X(K) \subseteq \oo{\widetilde{X}^*}$ denote the set  of all geodesics with respect to $X$ representing elements in $K$.  In a hyperbolic group, the language of all geodesics, $Geo_X(G)$, is rational for every generating set $X$. 
We say that a word $u\in \widetilde B^*$ is in \emph{normal form} if it is of the form $vb_i$, for some freely reduced word $v\in \widetilde A^*$ and $i\in [m]$. Clearly, for every $u\in \widetilde B^*$, there is a unique $\overline u\in \widetilde B^*$ in normal form such that $\overline u\pi=u\pi$. Notice that, when the word $u$ belongs to $\widetilde A^*$ this corresponds to free reduction.

Given two words $u, v \in \widetilde{X}^*$ we write $u \equiv v$ to emphasize that $u$ and $v$ are equal as words, while $u\pi = v\pi$ will be written to mean that they represent the same group element. 
We write $u \doteq u_1\ldots u_n$ if $u_1 \equiv u_1\ldots u_n$ with $u_1,\ldots, u_n \in \widetilde{X}$. For all $1 \leq i \leq j \leq n$, we write then $u^{[i,j]} = u_iu_{i+1}\ldots u_j$ and $u^{[j]} = u^{[1,j]}$.
Given a language $L$, we denote by $\Cyc(L)$ the language of all cyclic permutations of words in $L$. If $L$ is rational (resp. context-free), then $\Cyc(L)$ is also rational (resp. context-free \cite[Exercise 6.4 c)]{[HU79]}).

In \cite[Proposition 3.1]{[HR02]}, it is proved that if $u$ and $v$ are words in a $\delta$-hyperbolic group with $u\pi= v\pi$, $u$ is geodesic and $v$ is $(\lambda, \varepsilon)$-quasigeodesic, then $u$ and $v$ \emph{boundedly asynchronously $K$-fellow travel for some constant $K$ and some asynchronicity bound $M$}, where $K$ and $M$ depend only on $\lambda$, $\varepsilon$ and $\delta$.  With our notation, it follows from their proof that, given $\lambda, \varepsilon$, there exists a $K$ such that for all geodesic words $u$ and all $(\lambda,\varepsilon)$-quasigeodesic $v$ such that $u\pi=v\pi$, 
there is a (not necessarily strictly) increasing function $h:\{0,\ldots, |v|\}\to \{0,\ldots, |u|\}$ such that $h(0) = 0$, $h(|v|) = |u|$ and
\begin{align}\label{holtrees}
d(v^{[i]}\pi,u^{[h(i)]}\pi)\leq K \quad \text{ and } \quad |h(i)-h(i-1)|\leq 2K+1\end{align}
for $i \in [|v|]$.
We will denote the boundedly asynchronously fellow travel constant by $K(\lambda, \varepsilon, \delta)$ throughout the paper.
In particular, $(\lambda, \varepsilon)$-quasigeodesics and geodesics representing the same elements are at  Hausdorff distance at most $K(\lambda, \varepsilon,\delta)$. 

For a finite alphabet $A$, we say that $\mathfrak T=(Q, q_0,F, \delta, \lambda)$ is a \emph{finite state $A$-transducer} if $Q$ is a finite set of states, $q_0\in Q$ is the initial state, $F\subseteq Q$ is a set of final states,
$\delta:Q\times A\to Q$ and $\lambda:Q\times A\to A^*$ are mappings. We will write 
$x\xrightarrow{c|d} y$ to mean that $(x,c)\delta=y$ and $(x,c)\lambda=d$.
Given $L \subseteq A^*$, we write 
$$\mathfrak T(L) = \{ w_1\ldots w_n \mid \exists\, \mbox{a path\;} q_0\xrightarrow{a_1|w_1} q_1\xrightarrow{a_2|w_2} \cdots \xrightarrow{a_{n}|w_{n}} q_n \in T\mbox{ with }a_1\ldots a_n \in L\}.$$

\begin{theorem}\label{transducer}
Let $G$ be a f.g. virtually free group, X be a generating set and $L$ be a rational language of $(\lambda,\varepsilon)$-quasigeodesic words over $\tilde X$, for some (fixed) values $\lambda$ and $\varepsilon$. Then $Geo_X(L\pi)$ is  an (effectively computable) rational language.
\end{theorem}
\noindent\textit{Proof.} 
Let $K$ be the constant from (\ref{holtrees}) and $Q$ be the set of all geodesic words over $\tilde X$ of length at most $K$. Consider the finite transducer $\mathfrak T$
with set of vertices $Q$, edges $w\xrightarrow{c|u} v$ for $c\in \tilde X$, $u\in \tilde X^*$ a geodesic word of length at most $2K+1$, and $v$ a geodesic word representing $(u^{-1}wc)\pi$, and with the empty word being the initial and (unique) final state.
 We claim that $(\mathfrak T(L))\pi\subseteq L\pi$ and that $Geo_X(L\pi)\subseteq \mathfrak T(L)$, and so $Geo_X(L\pi)=Geo_X(G)\cap \mathfrak T(L)$ is a rational language.

Let $u \in \mathfrak T(L)$. There must be some $v \doteq v_1\ldots v_n\in L$ and a path of the form
$$\varepsilon = p_0 \xrightarrow{v_1|u_1} p_1\xrightarrow{v_2|u_2} p_2\cdots \xrightarrow{v_{n-1}|u_{n-1}} p_{n-1} \xrightarrow{v_{n}|u_{n}} p_n = \varepsilon$$
in $\mathfrak T$ with $u \equiv u_1\cdots u_n$. But then $(u_i^{-1}p_{i-1}v_i) = p_i$ for $i \in [n]$ and it follows easily by induction that 
$1= p_n\pi = (u_{n}^{-1}\cdots u_1^{-1}v_1\cdots v_n)\pi$, i.e., $u\pi=v\pi\in L\pi$. Therefore $(\mathfrak T(L))\pi\subseteq L\pi$.

Now, let $u \doteq u_1\ldots u_k\in Geo_X(L\pi)$. Then $u\pi = v\pi$ for some quasigeodesic $v \doteq v_1\ldots v_n \in L$.
Let $h:[n]\to [k]$ be the function from (\ref{holtrees}). For $i = 0,\ldots,n$, let $w_i \in Geo_X(((u^{[h(i)]})^{-1}v^{[i]})\pi)$.
We claim that there is a path in $\mathfrak T$ of the form
$$\varepsilon = w_0 \xrightarrow{v_1|u^{[h(1)]}} w_1\xrightarrow{v_2|u^{[h(1)+1,h(2)]}}\cdots \xrightarrow{v_{n}|u^{[h(n-1)+1,h(n)]}} w_n = \varepsilon.$$

Indeed, it follows from (\ref{holtrees}) that $w_i \in Q$ and $|u^{[h(i-1)+1,h(i)]}| \leq 2K+1$ for $i \in [n]$. The edges are well defined since
$$\begin{array}{lll}
((u^{[h(i-1)+1,h(i)]})^{-1}w_{i-1}v_i)\pi&=&((u^{[h(i-1)+1,h(i)]})^{-1}(u^{[h(i-1)]})^{-1}v^{[i-1]}v_i)\pi\\
&&\\
&=&((u^{[h(i)]})^{-1}v^{[i]})\pi = w_i\pi
\end{array}$$
holds for $i \in [n]$.

Hence 
$$u = u_1\ldots u_k = u_1\ldots u_{h(n)} \in \mathfrak T(v_1\ldots v_n) = \mathfrak T(v) \subseteq \mathfrak T(L)$$
and so $Geo_X(L\pi)\subseteq \mathfrak T(L)$. Therefore $Geo_X(L\pi)=Geo_X(G)\cap \mathfrak T(L)$ is a rational language.
\qed\\

\begin{remark}\label{rmk transduction}
The theorem above is stated in terms of rational languages but works in the exact same way for any class of languages preserved by rational transduction, such as the class of context-free languages.
\end{remark}

We will now prove that  the language of normal forms of words consists of quasigeodesics. For $i\in [m]$ we will denote by $\varphi_i$ the automorphism of $F$ defined by  $u\varphi_i = b_iub_i^{-1}$.

\begin{lemma}\label{lengthC}
Let $w\in \widetilde B^*$, $M=\max\{|a\varphi_i|_A\mid a\in A, i\in [m]\}$, $N=\max\{|u|_A: \exists i,j,k\in [m] : b_ib_j=ub_k\}$ and $C=\max\{M,N\}$. 
Then, $|\overline w| \leq C|w|$.
\end{lemma}
\noindent\textit{Proof.} Let $w\in  \widetilde B^*$. We proceed by induction on $|w|$. If $|w|=0$, then $\overline w=w$ and so $|\overline w|=|w|$. Now assume that the result holds for all  words of length up to some $n$ and let $w\in \widetilde B^*$ be such that $|w|=n+1$. Then $w=ux$ for some $x\in \widetilde B$ and we may write  $\overline u=vb_j$. From the induction hypothesis, it follows that $|vb_j|\leq C|u|$. If $x\in \widetilde A$, we have that $w\pi=(\overline ux)\pi= (vb_jx)\pi=(v\pi) (x\pi\varphi_j)b_j$, 
hence $\oo{w} = \oo{v(x\varphi_j)}b_j$
 and $$|\overline w|\leq |v| + |\overline{x\varphi_j}| +1 =|vb_j| + M\leq C|u|+C=C(|u|+1)=C|w|.$$
If $x\in \widetilde{\{b_1,\ldots b_m\}}$, say $x=b_r$,  then $\overline{b_jx}=yb_s$ for some $s\in [m]$ and $y\in F$ such that $|y|\leq N$.
Hence
$w\pi=(\overline ux)\pi= (vb_jx)\pi=(vyb_s)\pi$, yielding $\overline w= \overline{vy}b_s$ and $$|\overline w|\leq |v| + |y| +1 \leq |vb_j| + N\leq C|u|+C=C(|u|+1)=C|w|.$$
\qed\\

\begin{corollary} \label{qg}
Every word in normal form is a $(C,0)$-quasigeodesic for $C$ defined as in the preceding lemma.
\end{corollary}
\noindent\textit{Proof.} Let $w\in \widetilde B^*$ be a word in normal form. We have to show that any subword of $w$ of length $k$ has geodesic length at least $\frac k C$. Since any subword of $w$ is a word in normal form, we only need to prove that a word $u$ in normal formal has geodesic length of at least $\frac {|u|}C$. This follows from Lemma \ref{lengthC}, since, for a word $u$ in normal form, letting $v$ be a geodesic word such that $v\pi=u\pi$, we have that $\overline v=u$, and so $|u|\leq C|v|$, i.e. $|v|\geq \frac{|u|}{C}$.
\qed\\

\begin{lemma}\label{generators}
Let $G$ be a hyperbolic group, $X,Y$ be two generating sets and $\pi_X:\tilde X\to G$ and $\pi_Y:\tilde Y\to G$ be the natural surjective homomorphisms and put $N_{X,Y}=\max\{d_Y(1,x)\mid x\in X\}$. If  $u \doteq x_1 \cdots x_n$ is a geodesic word in $\Gamma_X(G)$, then a word of the form $v=v_1\cdots v_n$, where $v_i$ is a geodesic word in $\Gamma_Y(G)$ representing $x_i\pi$, is a $(N_{X,Y}^2,2N_{X,Y}^3)$-quasigeodesic in $\Gamma_Y(G)$.
\end{lemma}
\noindent\textit{Proof.} We have to prove that, for all $1 \leq i \leq j \leq |v|,$ $$j-i\leq N_{X,Y}^2d_Y(v^{[i]}\pi_Y,v^{[j]}\pi_Y)+2N_{X,Y}^3.$$
$$
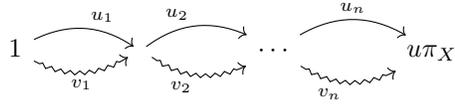
\begin{tikzcd}[sep=large]
   1\ar[r,bend left,"u_1"]\ar[r,swap,bend right, squiggly,"v_1"]  &\ar[r,bend left,"u_2"]\ar[r,swap,bend right, squiggly,"v_2"]  &\cdots   \ar[r,bend left,"u_n"]\ar[r,swap,bend right, squiggly,"v_n"] &u\pi_X
\end{tikzcd}
$$
Let $1 \leq i \leq j \leq |v|$. Define $k_i$ to be the largest integer such that $v_1\cdots v_{k_i}$ is a prefix of $v^{[i]}$ and $k_j$ to be the smallest integer 
such that 
 $v^{[j]}$ is a prefix of 
 $v_1\cdots v_{k_j}$. Notice that, for all $i$, $|v_i|\leq N_{X,Y}$.
Then, we have that:
\begin{align*}
j-i&\leq |v_{k_i+1}\cdots v_{k_j}|\\
&\leq |k_i-k_j|N_{X,Y}\\
&= N_{X,Y}d_X((v_1\cdots v_{k_i})\pi_Y,(v_1\cdots v_{k_j})\pi_Y)\\
&\leq N_{X,Y}^2d_Y((v_1\cdots v_{k_i})\pi_Y,(v_1\cdots v_{k_j})\pi_Y) \\
&\leq N_{X,Y}^2\left(d_Y((v_1\cdots v_{k_i})\pi_Y,v^{[i]}\pi_Y) + d_Y(v^{[i]}\pi_Y,v^{[j]}\pi_Y)+ d_Y(v^{[j]}\pi_Y,(v_1\cdots v_{k_j})\pi_Y)\right)\\
&\leq N_{X,Y}^2\left(N_{X,Y} + d_Y(v^{[i]}\pi_Y,v^{[j]}\pi_Y)+ N_{X,Y})\right)\\
&=  N_{X,Y}^2d_Y(v^{[i]}\pi_Y,v^{[j]}\pi_Y)+2N_{X,Y}^3.
\end{align*}
\qed\\

We can now combine the previous results to prove a generalization of Benois's Theorem for virtually free groups.

\begin{corollary}\label{benois vfree}
Let $G$ be a f.g. virtually free group and $K\subseteq G$. The following are equivalent:
\begin{enumerate}
\item $K\in \Rat(G)$.
\item $Geo_X(K)$ is a rational language for some finite generating set $X$.
\item $Geo_X(K)$ is a rational language for every finite generating set $X$.
\end{enumerate} 
Moreover, the constructions are effective.
\end{corollary}
\noindent\textit{Proof.} It is clear from the definitions that $3\implies 2 \implies 1$. We will prove that $1\implies 2$ and that $2\implies 3$. 

Compute the (rational) language $L$ of normal forms of $K$: this can be done by computing rational subsets $L_i$ of $F$ such that 
$K=\bigcup_{i\in [m]}L_ib_i$ (see \cite[Proposition 4.1]{[Sil02b]})
and then using Benois's theorem to compute the language of reduced words $\overline{L_i}$ representing elements in $L_i$. We then obtain that $$L=\bigcup_{i\in [m]} \overline{L_i}b_i.$$
In view of Corollary \ref{qg},  the language $L$ of normal forms representing elements in $K$ is a language of $(C,0)$-quasigeodesics over $\tilde B$ such that $L\pi=K$. By Theorem
\ref{transducer}, $Geo_B(K)$ is rational, so we have that $1\implies 2$.

Now, using Lemma \ref{generators}, we have that, given two finite generating sets $X,Y$ and replacing every edge of an automaton representing $Geo_X(K)$ by a path labeling a geodesic word over $\tilde Y$ representing the letter from $X$ labelling the edge, the language recognized by the new automaton will be a language $L$ of $(N_{X,Y}^2,2N_{X,Y}^3)$-quasigeodesic words over $\tilde Y$ such that $L\pi_Y=K$. Hence, $Geo_Y(K)$ is rational by Theorem \ref{transducer}.
\qed\\

\begin{remark} Similarly to what happens in Remark \ref{rmk transduction}, the equivalence between $2$ and $3$ holds for any class of languages closed under rational transductions.
Since rationality (resp. context-freeness) of the language of geodesics representing a given subset is independent of the generating set we will usually say that, for a subset $K$, $Geo(K)$ is rational (resp. context-free) to mean that $Geo_X(K)$ is rational (resp. context-free) for some (and so, for every) finite generating set $X$.
\end{remark}

We define $w$ to be a \emph{fully $(\lambda,\varepsilon)$-quasireduced word} if $w$ and all of its cyclic conjugates are $(\lambda,\
\varepsilon)$-quasigeodesic words.

We now present three  results from \cite{[HRR11]}  and \cite{[Car22c]}:

\begin{lemma}\cite[Lemma 16]{[HRR11]}\label{lema holt qr}
If $u$ and $v$ are fully $(\lambda,\varepsilon)$-quasireduced words representing conjugate elements of a $\delta$-hyperbolic group, then either $\max(|u|,|v|)\leq \lambda(8\delta+2K+\varepsilon+1)$ or there exist cyclic conjugates $u'$ and $v'$ of $u$ and $v$ and a word $\alpha$ with $(\alpha u'\alpha^{-1})\pi=v'\pi$ and $|\alpha|\leq 2(\delta+K)$, where $K$ is the boundedly asynchronous fellow travel constant satisfied by $(\lambda,\varepsilon)$-quasigeodesics with respect to geodesics. 
\end{lemma}
\begin{proposition}\cite[Proposition 18]{[HRR11]}\label{prop holt geo}
Let $u$ be a geodesic word in a $\delta$-hyperbolic group $G$ with $\delta \geq 1$. Then we have that $u\equiv u_1u_2u_3$, where $(u_3u_1)\pi=\alpha\pi$ for some word $\alpha$ with $|\alpha|\leq \delta$, and $u_2\alpha$ is fully $(1, 3\delta+1)$-quasireduced.

In other words, the word $u'\equiv u_1u_2\alpha\alpha^{-1}u_3$ obtained by insertion of $\alpha\alpha^{-1}$ into $u$ can be split as $u_1'u_2'u_3'$ such that $(u_3'u_1')\pi=1$ and $u_2'=u_2\alpha$ is fully $(1,3\delta+1)$-quasireduced.
\end{proposition}

Let $G$ be a hyperbolic group with generating set $A$. 
Given $g,h,p\in G$, we define the {Gromov product of $g$ and $h$} taking $p$ as basepoint by
$$(g|h)_p^A=\frac 12 (d_A(p,g)+d_A(p,h)-d_A(g,h)).$$ We will often write $(g|h)_p$ to denote $(g|h)_p^A$, when the generating set is clear from context.

\begin{lemma}\cite[Lemma 4.1]{[Car22c]}\label{meu lema}
Let $H$ be a hyperbolic group, $u,v\in H$ and $p\in \N$. Then the following are equivalent:
\begin{enumerate}[(i)]
\item $(u|v)_1\leq p$
\item for any geodesics $\alpha$ and $\beta$ from $1$ to $u^{-1}$ and $v$, respectively, we have that the concatenation
$$\xymatrix{
1 \ar[rr]^{\alpha} && u^{-1} \ar[rr]^{\beta} && u^{-1}v
}$$
is a $(1,2p)$-quasi-geodesic
\item there are geodesics $\alpha$ and $\beta$ from $1$ to $u^{-1}$ and $v$, respectively, such that the concatenation
$$\xymatrix{
1 \ar[rr]^{\alpha} && u^{-1} \ar[rr]^{\beta} && u^{-1}v
}$$
is a $(1,2p)$-quasi-geodesic
\end{enumerate}
\end{lemma}

\begin{lemma}\label{lemma geoquasireduced}
Let $G$ be a $\delta$-hyperbolic group and $g\in G$ be an element having a fully $(1,\varepsilon)$-quasireduced representative word $v$. Then, all geodesic words $w$ such that $w\pi=g$ are fully  $(1,\varepsilon +2K(1, \varepsilon,\delta)+2)$-quasireduced words.
\end{lemma}
\noindent\textit{Proof.} Put $K=K(1, \varepsilon,\delta)$.
Any geodesic $w$ is clearly a $(1,\varepsilon +2K+2)$-quasigeodesic. Now let $w=w_1w_2$ and consider the cyclic permutation $w'=w_2w_1$ of $w$. We have to prove that $w'$ is a   $(1,\varepsilon +2K+2)$-quasigeodesic. Consider the bigon with sides $w=w_1w_2$ (top side) and $v$ (bottom side). Since $v$ and $w$ are at Hausdorff distance at most $K$, then there is a vertex on the bottom side at a distance at most $K+1$ from the vertex reached after reading $w_1$ on the top side (the $+1$ comes from the possibility that the closest point of the bottom side to the vertex on the top side might not be a vertex itself) and so there is a geodesic word $\alpha$ of length at most $K+1$ and words $v_1,v_2$ such that $v\equiv v_1v_2$, $v_1\pi=(w_1\alpha)\pi$ and $v_2\pi=(\alpha^{-1}w_2)\pi$.

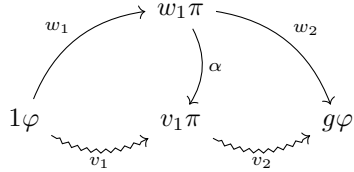
$$
\begin{tikzcd}[sep=large]
&w_1\pi \ar[dr,bend left,"w_2"]\ar[d,bend left,"\alpha"]&\\
   1\varphi \ar[r,swap,bend right, squiggly,"v_1"] \ar[ur,bend left,"w_1"] &v_1\pi \ar[r,swap,bend right, squiggly,"v_2"] & g\varphi
\end{tikzcd}
$$
We have that $(w_2w_1)\pi=(\alpha v_2v_1\alpha^{-1})\pi$ and so 
\begin{align*}
d(1,(v_2v_1)\pi)
&=d(1,(\alpha^{-1} w_2w_1\alpha)\pi)\\
&\leq 2|\alpha|+d(1,(w_2w_1)\pi)\\
&\leq d(w_2^{-1}\pi,w_1\pi)+2(K+1).\numberthis \label{eqn}
\end{align*}

Hence, using (\ref{eqn}) and the facts that $w\equiv w_1w_2$ is geodesic, $v\equiv v_1v_2$ is a $(1,\varepsilon)$-quasigeodesic (and so $|v|\leq d(1,(v_1v_2)\pi)+\varepsilon$), and that $v_2v_1$ is a $(1,\varepsilon)$-quasigeodesic (and so $|v_2v_1|\leq d(1,(v_2v_1)\pi)+\varepsilon$), we have that 
\begin{align*}
(w_2^{-1}\pi|w_1\pi)_1&=\frac 12(d(1,w_2^{-1}\pi)+d(1,w_1\pi)-d(w_2^{-1}\pi,w_1\pi))\\
&=\frac 12(|w_2|+|w_1|- d(w_2^{-1}\pi,w_1\pi))\\
&\leq\frac 12(|w|-d(1,(v_2v_1)\pi)+2(K+1))\\
&\leq \frac 12(|v|-d(1,(v_2v_1)\pi)+2(K+1))\\
&=\frac 12(|v_2v_1|-d(1,(v_2v_1)\pi)+2(K+1))\\
&\leq \frac \varepsilon 2+(K+1)\\
\end{align*}
From Lemma \ref{meu lema}, it follows that $w_2w_1$ is a $(1,\varepsilon +2K(1, \varepsilon,\delta)+2)$-quasigeodesic.

\qed\\

For convenience, we will denote $K(1, 3\delta+1,\delta)$ by $R$: this should cause no confusion, as the group, and so $\delta$, will be fixed. Recall that, for a subset $K\subseteq G$, we denote the set of all conjugates of elements of $K$ by $\alpha(K)$.
\begin{proposition}\label{uma qr}
Let $G$ be a virtually free group and $K\in \Rat(G)$. There is an effectively constructible rational language $L_K$ such that $L_K\pi\subseteq \alpha(K)$ and, for every element $g\in K$, there is at least one fully $(1,3\delta +2R+3)$-quasireduced word in $L_K$
 representing a conjugate of $g$.
\end{proposition}
\noindent\textit{Proof.}  Since $K$ is necessarily contained in some finitely generated subgroup of $G$, we may assume that $G$ is finitely generated. By Corollary \ref{benois vfree}, we can construct a finite state automaton recognizing $Geo_B(K)$, where $B$ is our standard generating set for $G$. Let $\delta$ be a hyperbolicity constant for $G$, $L_K=Geo_B(\Cyc(Geo_B(K))\pi)$ and $S=\{\alpha_1,\ldots, \alpha_n\}$  be the set of all words in $\tilde B^*$ of length at most $\delta$. We claim that $L_K$ has the desired properties.

The language $Geo_B(K)$ is rational in view of Corollary \ref{benois vfree}, and so $\Cyc(Geo_B(K))$ is rational. Hence, $\Cyc(Geo_B(K))\pi$ is a rational subset and $Geo_B(\Cyc(Geo_B(K))\pi)$ is rational by Corollary \ref{benois vfree}, which also implies that the construction is effective. 
 Also, $L_K\pi\subseteq \alpha(K)$, since, for every word $v\in L_K$, there is a word  $u\in \Cyc(Geo_B(K))$ such that $v\pi=u\pi$ and every word in $\Cyc(Geo_B(K))$  represents a conjugate of an element in $K$.

Now, let $g\in K$ and $u\in Geo_B(K)$ be a geodesic such that $u\pi=g$. Then, by Proposition \ref{prop holt geo}, there is some $i\in [n]$ such that $u\equiv u_1u_2u_3$, where $(u_3u_1)\pi=\alpha_i\pi$ , and $u_2\alpha_i$ is fully $(1, 3\delta+1)$-quasireduced. But, $(u_2\alpha_i)\pi=(u_2u_3u_1)\pi$ and $u_2u_3u_1\in \Cyc(Geo_B(K))$. Now,  any geodesic word representing $(u_2\alpha_i)\pi=(u_2u_3u_1)\pi$ belongs to $Geo_B(\Cyc(Geo_B(K))\pi)$ and,  by Lemma \ref{lemma geoquasireduced}, it is a fully $(1,3\delta +2R+3)$-quasireduced word. 
 \qed\\

\begin{theorem}\label{thmqr}
Let  $G$ be a finitely generated virtually free group and $K\in \Rat(G)$. There exists a context-free language $L'$ such that $L'\pi\subseteq \alpha(K)$ and $L'$ contains all the fully $(1,3\delta +2R+3)$-quasireduced words representing elements in $\alpha(K)$.\end{theorem}
\noindent\textit{Proof.} 
Let $L_K$ be the language from Proposition \ref{uma qr}, $$S=\{g\in G\mid d_B(1,g)\leq 2\delta+2K(1,3\delta+2R+3,\delta)\}$$ (notice that $S$ is finite), and $L=\Cyc(L_K)$. Since $L_K\pi\subseteq \alpha(K)$, then $L\pi\subseteq\alpha(K)$.

By brute force, we build the set $Q$ of all fully $(1,3\delta +2R+3)$-quasireduced words of length at most $ 11\delta +2K(1,3\delta+2R+3,\delta)+2R+4$ representing an element of $\alpha(K)$: it can be checked whether a word is a  quasigeodesic, and so it can be checked whether a word is fully quasireduced or not and, in case it is, we check if it belongs to $\alpha(K)$ using the main result from \cite{[LS11]}.

 Fix some $\beta\in S$. Since $L\pi\in \Rat(G)$, then $L\pi\in \CF(G)$ (see  \cite[Lemma 4.2]{[Her91]}) and 
 $$S_\beta:=\beta (L\pi)\beta^{-1}=
L\pi\lambda_\beta\in \CF(G),$$
which follows from Lemma \ref{square} by taking $M=M'=G$, $\tau=\lambda_\beta^{-1}$ and $T=L\pi$ and the fact that context-free languages are closed under inverse morphism.

Put
$L_\beta= S_\beta\pi^{-1}$. All words from $L_\beta$ represent a conjugate (by $\beta$) of an element in $L\pi$, and so all words in $L_\beta$ represent an element conjugate to an element in $K$, i.e.,  $L_\beta\pi\subseteq \alpha(K)$.

We claim that the language  $$L'=\bigcup_{\beta\in S}\Cyc(L_\beta) \cup Q$$ has the desired properties.
Clearly, it
is context-free and  $L'\pi\subseteq\alpha(K)$. We claim that it contains all the fully 
$(1,3\delta +2R+3)$-quasireduced words representing an element in $\alpha(K)$.

Let $v$ be a  fully $(1,3\delta +2R+3)$-quasireduced word representing an element in $\alpha(K)$.  We know that there is at least one fully $(1,3\delta +2R+3)$-quasireduced word $u\in L_K$ such that $u\pi\sim v\pi$ by Proposition \ref{uma qr}.

From Lemma \ref{lema holt qr}, it follows that either $\max(|u|,|v|)\leq   11\delta +2K(1,3\delta+2R+3,\delta)+2R+4$ or  there exist cyclic conjugates $u'$ and $v'$ of $u$ and $v$ and a word $\beta$ with $(\beta u'\beta^{-1})\pi=v'\pi$ and $|\beta|\leq 2\delta+2K(1,3\delta+2R+3,\delta)$. In the first case, we have that $v\in Q$, and so, $v\in L'$. So, assume that $|v|> 11\delta +2K(1,3\delta+2R+3,\delta)+2R+4$ and that there exist some $\beta\in S$ and cyclic permutations $u'$ and $v'$ of $u$ and $v$  with $\beta (u'\pi)\beta^{-1}=v'\pi$. In this case $u'\in L$ and $\beta (u'\pi)\beta\in S_\beta$, thus $v'\in L_\beta$ and $v\in \Cyc(L_\beta)$.

\qed

\begin{theorem} \label{alpha cf}
Let $G$ be a finitely generated virtually free group and $K\in \Rat(G)$. Then $Geo(\alpha(K))$ is context-free. \end{theorem}
\noindent\textit{Proof.} We will show that $Geo_B(\alpha(K))$ is context-free. Let $\delta$ be the maximum between $1$ and the hyperbolicity constant of $G$ (so $G$ is $\delta$-hyperbolic and $\delta \geq1$). It suffices to prove that there exists a context-free language $L$ such that $L\pi\subseteq \alpha(K)$ and 
 $Geo_B(\alpha(K))\subseteq L$, since, in that case $Geo_B(\alpha(K))=L\cap Geo_B(G)$ and context-free languages are closed under intersection with rational languages.
 
 Let $L'$ be the language given by Theorem \ref{thmqr}. For every $\beta\in \widetilde A^*$, the language $L'\cap\widetilde B^*\beta$ is context-free and then, so is the language $L_\beta''$ obtained by removing $\beta$ from the end of every word in $L'\cap\widetilde B^*\beta$.
 By the Muller-Schupp Theorem, $\{1\}\in \CF(G)$, and so $\{\beta\pi\}\in \CF(G)$, by Lemma \ref{lemma herbst}.
  Hence,
 the language $\beta\pi\pi^{-1}\subseteq \widetilde B^*$ is context-free and so is the language $\beta\pi\pi^{-1}\#\subseteq (\widetilde B\cup\#)^*$. Moreover,   $$L_2=\{u_1\#u_3\mid u_3u_1\in\beta\pi\pi^{-1}\}=\Cyc(\beta\pi\pi^{-1}\#)\subseteq  (\widetilde B\cup\#)^*$$
 is context-free. 
 Since context-free languages are closed under substitution, the language $$L_\beta=\{u_1u_2 u_3 \mid u_2\in L_\beta'', (u_3u_1)\pi=\beta\pi\}$$ obtained by replacing the symbol $\#$ by $L_\beta''$ in $L_2$ is context-free.
 We claim that the language $$L=\bigcup_{|\beta|\leq \delta} L_\beta$$ is context-free and  that $L\pi\subseteq \alpha(K)$ and 
 $Geo_B(\alpha(K))\subseteq L$. It is obvious that $L$ is context-free. Let $u_2\in L_\beta''$ and $u_1,u_3$ be such that $(u_3u_1)\pi=\beta\pi$ for some $\beta$ with $|\beta|\leq \delta$.
 Then $(u_1\pi)^{-1}(u_1u_2u_3)\pi(u_1\pi)=(u_2\beta)\pi$, and so $(u_1u_2u_3)\pi\sim (u_2\beta)\pi$. Since $u_2\in L_\beta''$, then $u_2\beta\in L'$ and $L'\pi\subseteq \alpha(K)$. Thus, $(u_1u_2u_3)\pi\in \alpha(K)$. Since $u_1u_2u_3$ is an arbitrary element of $L$, we have that $L\pi\subseteq \alpha(K).$ It remains to show that  $Geo_B(\alpha(K))\subseteq L$. Let $w\in Geo_B(\alpha(K))$. Then, by Proposition \ref{prop holt geo}, we have that  $w\equiv u_1u_2u_3$,  where $(u_3u_1)\pi=\beta\pi$ for some word $\beta$ such that $|\beta|\leq \delta$ and $u_2\beta$ is fully $(1,3\delta+1)$-quasireduced. It suffices to check that $u_2\in L_\beta''$, i.e., that $u_2\beta\in L'$. This follows from Theorem \ref{thmqr}, as every fully $(1,3\delta+1)$-quasireduced word is also fully $(1,3\delta +2R+3)$-quasireduced and $(u_2\alpha)\pi  = (u_2u_3u_1)\pi \sim w\pi \in \alpha(K)$. 
 \qed

\begin{corollary}\label{dgcp vfree}
Let $G$ be a virtually free group. Then the doubly generalized conjugacy problem is decidable.
\end{corollary}
\noindent\textit{Proof.} It amounts to deciding, on input $S,T\in \Rat(G)$, whether $Geo_B(\alpha(S))\cap Geo_B(T)=\emptyset$, which can be done since $Geo_B(\alpha(S))$ is context-free by Theorem \ref{alpha cf} and $Geo_B(T)$ is rational by Corollary \ref{benois vfree}.
\qed\\

\section*{Acknowledgements}
Both authors were partially supported by
CMUP, member of LASI, which is financed by national funds through FCT – Fundação
para a Ciência e a Tecnologia, I.P., under the projects with reference UIDB/00144/2020
and UIDP/00144/2020.

\bibliographystyle{plain}
\bibliography{Bibliografia}

 \end{document}